\documentclass[11pt]{amsart}
\usepackage{amssymb}

\newtheorem{lemma}{Lemma}[section]

\newtheorem{theorem}[lemma]{Theorem}

\newtheorem{proposition}[lemma]{Proposition}
\newtheorem{corollary}[lemma]{Corollary}
\newtheorem{definition}[lemma]{Definition}

\begin{document}
\title{L\^e-Greuel type formula for the Euler obstruction and applications}
\author{Nicolas Dutertre and Nivaldo G. Grulha Jr.}
\address{Universit\'e de Provence, Centre de Math\'ematiques et Informatique,
39 rue Joliot-Curie,
13453 Marseille Cedex 13, France.}
\email{dutertre@cmi.univ-mrs.fr}

\address{Universidade de S\~{a}o Paulo, Instituto de Ci\^{e}ncias Matem\'{a}ticas e de Computa\c{c}\~{a}o - USP  Av. Trabalhador S\~{a}o-carlense, 400 - Centro, Caixa Postal: 668 - CEP: 13560-970 - S\~{a}o Carlos - SP -
    Brazil.}
\email{njunior@icmc.usp.br}

\thanks{
N. Dutertre is supported by {\em Agence Nationale de la Recherche}
(reference ANR-08-JCJC-0118-01)\\
N. G. Grulha Jr. is supported by {\em Funda\c{c}\~{a}o de Amparo \`{a} Pesquisa do Estado de S\~{a}o Paulo }
(reference 2009/08774-0)}

\begin{abstract}The Euler obstruction of a function $f$ can be viewed as a generalization of the Milnor number for functions defined on singular spaces. In this work, using the Euler obstruction of a function, we give a version of the L\^{e}-Greuel formula for germs $f:(X,0) \to (\mathbb{C},0)$ and $g:(X,0) \to (\mathbb{C},0)$ of analytic functions with isolated singularity at the origin. Using this formula and results of Loeser, we also present an integral formula for the Euler obstruction of a function, generalizing a formula of Kennedy.
\end{abstract}

\maketitle
\markboth{N. Dutertre and N. G. Grulha Jr.}{L\^e-Greuel type formula for the Euler obstruction and applications }

\section{Introduction}
Let $f : (\mathbb{C}^n,0) \rightarrow (\mathbb{C},0)$ be an analytic function defined in a neighborhood of the origin. The Milnor fiber is the set $f^{-1}(\delta) \cap B_\varepsilon$, where $\delta$ is a regular value of $f$, $B_\varepsilon$ is the euclidian ball in $\mathbb{C}^n$ of radius $\varepsilon$ and where $0< \vert \delta \vert \ll \varepsilon \ll 1$. Many research works have been devoted to the study of the topology of the Milnor fiber. When $f$ has an isolated critical point at the origin, Milnor \cite{Mil} proved that $f^{-1}(\delta) \cap B_\varepsilon$ had the homotopy type of a wedge of spheres of dimension $n-1$. The number of spheres appearing in this wedge is called the Milnor number of $f$ and is denoted by $\mu(f)$. Milnor showed that $\mu(f)$ was equal to the topological degree of the map $\frac{\nabla f}{\Vert \nabla f \Vert} : S_\varepsilon \rightarrow S_1$, where $S_\varepsilon$ and $S_1$ are the euclidian spheres in $\mathbb{C}^n$ of radius $\varepsilon$ and $1$. Therefore $\mu (f)$ is equal to the dimension of the $\mathbb{C}$-algebra $\frac{\mathcal{O}_{\mathbb{C}^n,0}}{(\frac{\partial f}{\partial x_1},\ldots,\frac{\partial f}{\partial x_n})}$, where $\mathcal{O}_{\mathbb{C}^n,0}$ is the algebra of holomorphic germs at the origin, and also to the number of critical points of a Morsification of $f$.

Milnor's results were extended to the case of a complete intersection with isolated singularity $F=(f_1,\ldots,f_k) : (\mathbb{C}^n,0) \rightarrow (\mathbb{C}^k,0)$, $1<k<n$, by Hamm \cite{Ham}, who proved that the Milnor fiber
$F^{-1}(\delta) \cap B_\varepsilon$, $0< \vert \delta \vert \ll \varepsilon \ll 1$, had the homotopy type of a wedge of $\mu(F)$ spheres of dimension $n-k$. L\^e \cite{Le0} and Greuel \cite{Gre} proved the following formula:
$$\mu(F')+ \mu(F)= \hbox{dim}_{\mathbb{C}} \frac{\mathcal{O}_{\mathbb{C}^n,0}}{I}, $$
where $F' :(\mathbb{C}^n,0) \rightarrow (\mathbb{C}^{k-1},0)$ is the map with components $f_1,\ldots,f_{k-1}$ and
$I$ is the ideal generated by $f_1,\ldots,f_{k-1}$ and the $(k \times k)$-minors
$\frac{\partial(f_1,\ldots,f_k)}{\partial(x_{i_1},\ldots,x_{i_k})}$.
If we denote by $\Gamma_F$ the set of points where all these minors vanish then the above formula can be reformulated as follows:
$$\mu(F')+ \mu(F)=I_{\mathbb{C}^n,0} ( \{F'=0 \}, \Gamma_F),$$
where $I_{\mathbb{C}^n,0}(-,-)$ is the intersection multiplicity at the origin in $\mathbb{C}^n$.  When $F'=f : (\mathbb{C}^n,0) \rightarrow (\mathbb{C},0)$ and $F=(f,l)$, where $l$ is a generic linear function, the above formula is also called Teissier's lemma \cite{Te}. Note also that in this situation $\mu(f,l)$ is equal to the first Milnor-Teissier number $\mu'(f)$, namely the Milnor number of $f$ restricted to a generic hyperplane.

Using Teissier's lemma and tools from integral geometry, Langevin \cite{Lan} gave an integral formula for the sum $\mu(f)+ \mu'(f)$. In \cite{K}, Kennedy carried on Langevin's work and established an integral formula for the Milnor number $\mu(f)$ (see also \cite{Lo} for other integral formulas).

Our aim in this paper is to generalize some of the previous results to the singular case, namely  to replace $\mathbb{C}^n$ with an equidimensional reduced complex analytic space $X$. There exist several generalizations of the Milnor number of a function defined on a singular space. We will use the local Euler obstruction of a function defined in \cite{BMPS}.

Let us recall first that the local Euler obstruction at a point $p$ of an algebraic variety $X$,
denoted by ${\rm Eu}_{X}(p)$, was defined by MacPherson.
It is one of the main ingredients in his proof of Deligne-Grothendieck conjecture
concerning existence of characteristic classes for complex algebraic varieties
\cite{M}. An equivalent definition was given in  \cite{BS} by J.-P. Brasselet and M.-H.
Schwartz, using stratified vector fields.

The computation of local Euler obstruction is not so easy by
using the definition. Various authors propose formulae which make the
computation easier.  L\^{e} D.T.  and
B. Teissier provide a formula in terms of polar multiplicities \cite{LT1}.

In the paper  \cite{BLS}, J.-P. Brasselet, L\^{e} D. T. and J. Seade
give a Lefschetz type formula for the local Euler obstruction. The formula shows that
the local Euler obstruction, as a constructible function, satisfies the Euler condition
relatively to generic linear forms. A natural continuation of the result is the paper
by J.-P. Brasselet, D. Massey, A. J. Parameswaran and J. Seade \cite{BMPS},
whose aim is to understand what is the obstacle for the local Euler obstruction
to satisfy  the Euler condition relatively to analytic functions with isolated singularity
at the considered point. That is the role of the so-called local Euler obstruction
of $f$, denoted by ${\rm Eu}_{f,X}(0)$.

In this paper we introduce the number ${\rm B}_{f,X}(0)= {\rm Eu}_{X}(0) - {\rm Eu}_{f,X}(0)$.  When  $f$ is linear and generic, it gives ${\rm Eu}_{X}(0)$, hence B$_{f,X}(0)$ can be viewed as a generalization of the Euler obstruction, as it was first remarked by Jean-Paul Brasselet. Using the number B$_{f,X}(0)$, we give a version of the L\^{e}-Greuel formula for germs $f:(X,0) \to (\mathbb{C},0)$ and $g:(X,0) \to (\mathbb{C},0)$ of analytic functions with isolated singularity at the origin (Theorem \ref{IntNumbEuler}). Namely we prove:
$$\displaylines{
\qquad {\rm B}_{f,X}(0)-{\rm B}_{f,X^{g}}(0) = (-1)^d I_{X,0} (X^f, \overline{\Gamma^q_{f,g}}),
 \qquad \cr
}$$
where $X^f=X \cap f^{-1}(0)$, $X^g=X \cap g^{-1}(0)$ and $I_{X,0} (X^f, \overline{\Gamma^q_{f,g}})$ is the intersection multiplicity at the origin of $X^f$ and the following polar set:
$$\Gamma^q_{f,g} = \left\{ x \in V_i \ \vert \ \hbox{rank} (df_{\vert V_q}(x), dg_{\vert V_q}(x))< 2 \right\}.$$
Here $V_q=X_{\rm reg}$ is the set of smooth points of $X$.

When $g$ is a generic linear function, we obtain a singular version of Teissier's lemma (Theorem \ref{TeissierLemmaSing}).  As a corollary of this singular Teissier's lemma, we express ${\rm B}_{f,X}(0)$ in terms of relative polar multiplicities (Theorem \ref{EulerObstructionInterMult}), in the spirit of the formula of L\^e and Teissier mentioned above.

Then, using this last result and results of Loeser  \cite{Lo}, we also present an integral formula for  ${\rm B}_{f,X}(0)$, generalizing the formula of Kennedy (Corollary \ref{KennedyGen}) and a Gauss-Bonnet formula for the Milnor fiber of $f$ (Corollary \ref{GaussBonnetGen}).

The paper is organized as follows : in Section 2, we recall some facts about the Euler obstruction, the Euler obstruction of a function and the Euler obstruction of a $1$-form. In Section 3, we give the definition of the complex link, the definition of the radial index of a $1$-form and a useful relation due to Ebeling and Gusein-Zade between these two notions and the Euler obstruction of a $1$-form (Theorem \ref{Th4}). Section 4 is devoted to the proof of an auxiliary lemma in subanalytic geometry. In Section 5, we present two versions of the L\^e-Greuel formula (Theorem \ref{teo3.2} and Corollary \ref{IntNumbGen}) that we will use in the proof of our main results. We also give an application to $1$-parameter deformations (Corollary \ref{Deformation}). Section 6 contains the main results mentioned above.

Our result concern functions with isolated singularities, it would be interesting to try to generalize them to the non-isolated singularity case, using for instance technics of D. Massey \cite{Massey}.

The authors are grateful to M.A.S Ruas for suggesting them the application to $1$-parameter deformations, and to J.-P. Brasselet for his careful reading and for suggesting some improvements in this paper.

\section{The Euler obstruction}\label{qualquer}

The Euler obstruction was defined by MacPherson \cite{M} as a tool to prove the conjecture about existence and unicity of the Chern classes in the singular case. Since that the Euler obstruction was deeply investigated by many authors as Brasselet, Schwartz, Sebastiani, L\^{e}, Teissier, Sabbah, Dubson, Kato and others. For an overview about the Euler obstruction see \cite{B,BG}. Let us now introduce some objects in order to define the Euler obstruction.

For all this paper, let us consider the following setting. Let $(X,0) \subset (\mathbb{C}^N,0)$ be an equidimensional reduced complex analytic germ of dimension $d$ in an open set $U \subset \mathbb{C}^N$. We consider a complex analytic Whitney stratification $\{V_i\}$ of $U$ adapted to $X$ and we assume that $\{0\}$ is a stratum. We choose a small representative of $(X,0)$ such that $0$ belongs to the closure of all the strata. We will denote it by $X$ and we will write $X= \cup_{i=0}^q V_i$ where $V_0 = \{0\}$ and $V_q = X_{\rm reg}$, the set of smooth points of $X$. We will assume that the strata $V_0,\ldots,V_{q-1}$ are connected and that the analytic sets $\overline{V_0},\ldots,\overline{V_{q-1}}$ are reduced. We will set $d_i =\hbox{dim} V_i$ for $i \in \{1,\ldots,q\}$ (note that $d_q=d$).

Let $G(d,N)$ denote the Grassmanian of complex $d$-planes
in ${\mathbb C}^N$. On the regular part $X_{\rm reg}$ of $X$ the Gauss map
$\phi : X_{ \rm reg} \to U\times G(d,N)$ is well defined by $\phi(x) =
(x,T_x(X_{ \rm reg}))$.

\begin{definition}
The {\it Nash transformation} (or {\it Nash blow up}) $\widetilde X$ of $X$
is the closure of the image {\rm Im}$(\phi)$ in $ U\times G(d,N)$. It is a
(usually singular) complex analytic space endowed with an analytic
projection map $\nu : \widetilde   X \to X$
which is a biholomorphism away from $\nu^{-1}({{\rm Sing}}(X))\,.$
\end{definition}

The fiber of  the tautological bundle  ${\mathcal T}$ over $G(d,N)$, at the point $P \in G(d,N)$,
is the set of the vectors $v$ in the $d$-plane $P$. We still denote by ${\mathcal T}$ the
corresponding trivial extension bundle over $ U \times G(d,N)$. Let $\widetilde
T$ be the restriction of  ${\mathcal T}$ to $\widetilde   X$, with
projection map $\pi$.  The bundle $\widetilde   T$ on $\widetilde
X$ is called {\it the  Nash bundle} of $X$.

An element of $\widetilde   T$ is written $(x,P,v)$ where $x\in U$,
$P$ is a $d$-plane in ${\mathbb C}^N$ based at $x$ and $v$ is a vector in
$P$. We have the following diagram:
$$
\begin{matrix}
\widetilde   T & \hookrightarrow & {\mathcal T} \cr
{\pi} \downarrow & & \downarrow \cr
\widetilde   X & \hookrightarrow & U \times G(d,N) \cr
{\nu}\downarrow & & \downarrow \cr
X & \hookrightarrow & U. \cr
\end{matrix}
$$

Let us recall the original definition of the Euler obstruction, due to Mac\-Pherson \cite{M}.
Let $z=(z_1, \ldots,z_N)$ be local coordinates in ${\mathbb C}^N$ around $\{ 0 \}$, such that
$z_i(0)=0$. We denote by $ B_\varepsilon$ and $S_\varepsilon$ the ball and the sphere
centered at $\{ 0 \}$ and of radius $\varepsilon$ in ${\mathbb C}^N$. Let us consider the norm
$\Vert z \Vert = \sqrt{z_1 \overline z_1 + \cdots + z_N \overline z_N}$. Then the differential form
$\omega = d \Vert z \Vert ^2$ defines a section of the real vector bundle $T({\mathbb C}^N)^*$,
cotangent bundle on ${\mathbb C}^N$. Its pull back restricted to $ \widetilde   X$ becomes a section
denoted by $ \widetilde\omega$ of the dual bundle ${\widetilde   T}^*$. For $\varepsilon$ small enough,
the section  $ \widetilde\omega$ is nonzero over $\nu^{-1}(z)$ for $0 < \Vert z \Vert \le \varepsilon$.
The obstruction to extend $ \widetilde\omega$ as a nonzero section of ${\widetilde   T}^*$ from
$\nu^{-1}(S_\varepsilon)$ to $\nu^{-1}(B_\varepsilon)$, denoted by
$Obs({\widetilde   T}^*, \widetilde\omega)$ lies in $H^{2d}(\nu^{-1}( B_\varepsilon),
\nu^{-1}( S_\varepsilon);{\mathbb Z})$. Let us denote by ${\mathcal O}_{\nu^{-1}(B_\varepsilon),
\nu^{-1}( S_\varepsilon)}$ the orientation class in $H_{2d}(\nu^{-1}(B_\varepsilon),
\nu^{-1}(S_\varepsilon);{\mathbb Z})$.

\begin{definition}
The local Euler obstruction of $X$ at $0$ is the evaluation of $Obs({\widetilde   T}^*, \widetilde\omega)$ on
${\mathcal O}_{\nu^{-1}(\mathbb B_\varepsilon), \nu^{-1}(\mathbb S_\varepsilon)}$, i.e:
$${\rm Eu}_{X}(0) = \langle  Obs({\widetilde   T}^*, \widetilde\omega), {\mathcal O}_{\nu^{-1}(B_\varepsilon),
\nu^{-1}(S_\varepsilon)} \rangle.$$
\end{definition}

An equivalent definition of the Euler obstruction was given by Brasselet and Schwartz in the context of vector fields \cite {BS}.

We will need in this paper some results about the Euler obstruction
where this invariant is computed using hyperplane sections. The idea of studying the Euler obstruction using hyperplane sections appears in the works of Dubson
and Kato, but the approach we follow here is that of \cite {BLS,BMPS}.

\begin{theorem}[\cite{BLS}]\label{BLS}
Let $(X,0)$ and $\{V_{i}\}$ be given as before, then
for each generic linear form $l$, there is $\varepsilon_0$
such that for any $\varepsilon$ with $0<\varepsilon<
\varepsilon_0$ and $t_0\neq 0$ sufficiently small, the Euler obstruction of $(X,0)$ is equal to:
$${\rm Eu}_X(0)=\sum_{i =1}^q\chi \big(V_i\cap B_\varepsilon\cap l^{-1}(\delta) \big) \cdot
{\rm Eu}_{X}(V_i),$$ where $\chi$ denotes the Euler-Poincar\'e
characteristic, ${\rm Eu}_{X}(V_i)$ is the value of the Euler
obstruction of $X$ at any point of $V_i$, $i=1,\ldots,q$, and $0 < \vert \delta \vert \ll \varepsilon \ll 1$.
\end{theorem}

We define now an invariant introduced by Brasselet, Massey, Parameswa\-ran and Seade in \cite
{BMPS}, which measures in a way how far the equality given in
Theorem \ref{BLS} is from being true if we replace the generic
linear form $l$ with  some other function on $X$ with at most an
isolated stratified critical point at $0$. So let $f : X \rightarrow \mathbb{C}$ be  a holomorphic function which is the restriction of a holomorphic function $F : U \rightarrow \mathbb{C}$. A point $x$ in $X$ is a critical point of $f$ if it is a critical point of $F_{\vert V(x)}$, where $V(x)$ is the stratum containing $x$. We will assume that $f$ has an isolated singularity (or an isolated critical point) at $0$, i.e that $f$ has no critical point in a punctured neighborhood of $0$ in $X$.
In order to define the new invariant the authors constructed in \cite{BMPS} a stratified vector field on
$X$, denoted by $\overline{\nabla}_{X} f$. This vector field is homotopic
  to $\overline {\nabla} F\vert_X$ and  one has $\overline{\nabla}_{X} f (x)
\ne 0 $ unless $x = 0$.

Let $\tilde\zeta$ be the lifting of $\overline{\nabla}_{X} f$ as a
section of the Nash bundle $\widetilde   T$ over $\widetilde   X$
without singularity over $\nu^{-1}(X\cap S_\varepsilon)$. Let ${\mathcal O} (\tilde\zeta) \in H^{2n}\big(\nu^{-1}(X \cap
B_\varepsilon), \nu^{-1}(X \cap
S_\varepsilon)\big)$ be the obstruction cocycle
  to the extension of  $\widetilde  \zeta$
as a nowhere zero section of $\widetilde   T$ inside $\nu^{-1}
(X\cap B_\varepsilon)$.

\begin{definition}\label{demai}
The local Euler obstruction
$ {\rm Eu}_{f,X}(0)$ is the evaluation of ${\mathcal O}
(\tilde\zeta)$ on the fundamental class of the pair $(\nu^{-1}(X
\cap B_\varepsilon), \nu^{-1}(X \cap
S_\varepsilon))$.
\end{definition}

The following result compares
the Euler obstruction of the space $X$ with that of a function on
$X$ \cite {BMPS}.

\begin{theorem}\label{BMPS}
Let $(X,0)$ and $\{V_{i}\}$ given as before and let $ f : (X,0)  \to (\mathbb C,0)$ be a function with an isolated
singularity at $0$. For $0 < \vert \delta \vert \ll \varepsilon \ll 1$ we have:
$${\rm Eu}_{f,X}(0)={\rm Eu}_X(0) \,- \,\left(\sum_{i=1}^{q} \chi \big(V_i\cap B_\varepsilon\cap
f^{-1}(\delta) \big) \cdot  {\rm Eu}_X(V_i) \right).$$
\end{theorem}

Using the notation defined in the introduction, we have:

$${\rm B}_{f,X}(0)= \sum_{i=1}^{q} \chi \big(V_i\cap B_\varepsilon\cap
f^{-1}(\delta) \big) \cdot  {\rm Eu}_X(V_i).$$

In \cite{STV}, J. Seade et al. show that the Euler obstruction of $f$ is closely
related to the number of Morse points of a Morsification of $f$, as it is stated in the next proposition.

\begin{proposition}[\cite{STV}]\label{Prop-n-reg} Let $f:(X,0) \to (\mathbb{C},0)$ be
the germ of an analytic function with isolated singularity at the
origin.
Then: $${\rm Eu}_{f,X}(0)=(-1)^d n_{reg},$$where
$n_{reg}$ is the number of Morse points in $X_{\rm reg}$ in a stratified
morsification of $f$.
\end{proposition}

Let us consider the Nash bundle $\widetilde   T$ on $\widetilde X$.
The corresponding dual bundles of complex and real 1-forms  are
denoted, respectively, by $\widetilde T^* \buildrel {}\over \to \widetilde X $
and $\widetilde T^*_{\mathbb R} \buildrel {}\over \to \widetilde X $.

\begin{definition}\label{def2.1} Let $(X,0)$ and $\{V_\alpha\}$ as before.
Let $\omega$ be a (real or complex) 1-form on $X$,
 {\it i.e.}, a continuous section of either $T^*_{\mathbb R} \mathbb{C}^{N} \vert_X$ or $T^*\mathbb{C}^{N} \vert_X$.
{\it A singularity of  $\omega$ in the stratified sense}
means a point  $x$ where the kernel of
$\omega$  contains the tangent space of the corresponding stratum.
\end{definition}

This means that the pull back of the form to $V_\alpha$ vanishes at
$x$. Given a section $\eta$ of $T^*_{\mathbb R}{\mathbb C}^N |_A$, $A\subset V$,
there is a canonical way of constructing a
section $\tilde\eta$ of $\widetilde T^*_{\mathbb R}|_{\tilde A}$, $\tilde
A=\nu^{-1}A$, such that if $\eta$ has an isolated singularity at the point $0 \in X$ (in the
stratified sense), then we have a never-zero section $\widetilde
\eta$ of the dual Nash bundle $\widetilde T^*_{\mathbb R}$ over
$\nu^{-1}({{ S}_\varepsilon} \cap X) \subset \widetilde X$. Let $o(\eta) \in
H^{2d}(\nu^{-1}({{ B}_\varepsilon} \cap X), \nu^{-1}({{ S}_\varepsilon} \cap X); {\mathbb Z})$ be
the cohomology class of the obstruction cycle to extend this to a
section of  $\widetilde T^*_{\mathbb R}$ over $\nu^{-1}({{ B}_\varepsilon} \cap X)$.
Then we can define (c.f. \cite {EG1}):

\begin{definition} The {\it local Euler obstruction} of the real differential form $\eta$
at an isolated singularity is the integer ${\rm Eu}_{X,0}\ \eta$
obtained by evaluating the obstruction cohomology class $o(\eta)$
on the orientation fundamental cycle $[\nu^{-1}({{B}_\varepsilon} \cap X),
\nu^{-1}({{ S}_\varepsilon}\cap X)]$.
\end{definition}

MacPherson's local Euler obstruction $\hbox{Eu}_{X}(0)$
corresponds to taking  the differential $\omega= d{\Vert z \Vert}^{2}$ of the square of the  function
distance to $0$.

In the complex case, one can perform the same construction, using
the corresponding complex bundles. If $\omega$ is a complex
differential form, section of $T^*{\mathbb C}^N |_A$ with an isolated
singularity, one can define  the local Euler obstruction
$\hbox{Eu}_{X,0}\  \omega $.
Notice that, as explained in \cite{BSS} p.151, it is equal to the local Euler obstruction of its real part up to sign:
\begin{equation}\label{realpart}
{\rm Eu}_{X,0}\ \omega \,=\, (-1)^d {\rm Eu}_{X,0}{\rm Re}\ \omega .\notag
\end{equation}
This is an immediate consequence of the relation between the Chern
classes of a complex vector bundle and those of its dual. Remark also that when we consider the differential of a function $f$, we have the following equality (see \cite{EG1}):
$${\rm Eu}_{X,0}\ df=(-1)^d {\rm Eu}_{f,X}( 0).$$

We note that the idea to consider the (complex) dual Nash
bundle\index{Nash bundle}
 was already present in
\cite{Sa},  where Sabbah introduces a local
Euler obstruction ${\rm E\check u}_X (0)$ that satisfies ${\rm
E\check u}_X (0) = (-1)^{d}{\rm Eu}_X(0)$.  See also
 \cite[sec. 5.2]{Schu}.

\section{The complex link and the radial index}
In this section, we recall the definition of the complex link and of the radial index. We also present a formula of Ebeling and Gusein-Zade which expresses the radial index of a $1$-form  in terms of Euler characteristics of complex links and Euler obstructions.

The complex link is an important object in the study of the topology of complex analytic sets. It is analogous to the Milnor fibre and was studied first in \cite{Le1}. It plays a crucial role in complex stratified Morse theory (see \cite{GMP}) and appears in general bouquet theorems for the Milnor fibre of a function with isolated singularity (see \cite{Le2, Si, Ti}). It is related to the multiplicity of polar varieties and also the local Euler obstruction (see \cite{Db1, Db2, LT1, LT2}).  Let us recall briefly its definition. Let $M$ be a complex analytic manifold equipped with a Riemannian metric and let $Y \subset M$ be a complex analytic variety equipped with a Whitney stratification. Let $V$ be a stratum of $Y$ and let $p$ be a point in $V$. Let $N$ be a complex analytic submanifold of $M$ which meets $V$ transversally at the single point $p$. By choosing local coordinates on $N$, in some neighborhood of $p$ we can assume that $N$ is an Euclidian space $\mathbb{C}^k$.

\begin{definition}
The complex link of $V$ in $Y$ is the set denoted by ${\rm lk}^\mathbb{C}(V,Y)$ and defined as follows:
$$ {\rm lk}^\mathbb{C}(V,Y)= Y \cap N \cap B_\varepsilon \cap l^{-1}(\delta),$$
where $l : N \rightarrow \mathbb{C}$ is a generic linear form and $0< \vert \delta \vert \ll \varepsilon \ll 1$.
\end{definition}
The fact that the link of a stratum is well-defined, i.e independent of all the choices made to define it, is explained in \cite{LT2, Db2, GMP}. It is also independent of the embedding of the analytic variety $Y$ (see \cite{LT2}).  Now let $H$ be a smooth analytic hypersurface of $M$ transverse to $Y$. The intersections of $H$ with the stratum of $Y$ give a Whitney stratification of $Y \cap H$. We will need the following lemma.
\begin{lemma}\label{complexlink}
Let $V$ be a stratum of $Y$ such that $V \cap H \not= \emptyset$. We have ${\rm lk}^\mathbb{C}(V,Y)={\rm lk}^\mathbb{C}(V \cap H,Y \cap H)$.
\end{lemma}
{\it Proof.} Let $p$ be a point in $V \cap H$. Let $N \subset H$ be a complex analytic manifold of dimension ${\rm dim}\ H -{\rm dim} (V \cap H)$ that intersects $V \cap H$ transversally at $p$ in $H$. Since $H$ intersects $V$ transversally, we have:
$${\rm dim}\ H-{\rm dim}( V \cap H) = {\rm dim}\  M-1 -({\rm dim}\ V -1) = {\rm dim}\ M - {\rm dim}\ V.$$
 Hence $N$ intersects $V$ transversally at $p$ in $M$ as well. We conclude remarking that
 $ Y \cap N = Y \cap H \cap N$ because $N \subset H$. $\hfill \Box$

In \cite{EG1}, Ebeling and Gusein-Zade established relations between the local Euler obstruction of a $1$-form, its radial index and Euler characteristics of complex links. The radial index is a generalization to the singular case of the Poincar\'e-Hopf index. In order to define this index, let us consider first the real case.
Let $Z \subset \mathbb{R}^n$ be a closed subanalytic set equipped with a Whitney stratification $\{S_\alpha\}_{\alpha \in
\Lambda}$. Let $\omega$ be a continuous 1-form defined on $\mathbb{R}^n$. We say that a point $P$ in $Z$ is a zero (or a
singular point) of $\omega$
on $Z$ if it is a zero of $\omega_{\vert S}$, where $S$ is the stratum that contains $P$. In the sequel, we will define the radial
index of $\omega$ at $P$, when $P$ is an isolated zero of $\omega$ on $Z$. We can assume that $P=0$ and we denote by $S_0$ the stratum that contains $0$.

\begin{definition}
A 1-form $\omega$ is radial on $Z$ at $0$ if, for an arbitrary non-trivial subanalytic arc $\varphi : [0,\nu[
\rightarrow Z$ of class $C^1$, the value of the form $\omega$ on the tangent vector $\dot{\varphi}(t)$ is positive for $t$
small enough.
\end{definition}

Let $\varepsilon >0$ be small enough so that in the closed ball $B_\varepsilon$, the 1-form has no singular points on $Z \setminus \{0\}$. Let $S_0,\ldots,S_r$ be the strata that contain $0$
in their closure. Following Ebeling and Gusein-Zade, there exists a 1-form $\tilde{\omega}$ on
$\mathbb{R}^n$ such that :
\begin{enumerate}
\item The 1-form $\tilde{\omega}$ coincides with the 1-form $\omega$ on a neighborhood of $S_\varepsilon$.
\item The 1-form $\tilde{\omega}$ is radial on $Z$ at the origin.
\item In a neighborhood of each zero $Q \in Z \cap B_\varepsilon \setminus \{0\}$, $Q \in S_i$, ${\rm dim} \ S_i=k$, the
1-form $\tilde{\omega}$ looks as follows. There exists a local subanalytic diffeomorphism $h :(\mathbb{R}^n,\mathbb{R}^k,0)
\rightarrow (\mathbb{R}^n,S_i,Q)$ such that $h^* \tilde{\omega}=\pi_1^* \tilde{\omega}_1+\pi_2^*\tilde{\omega}_2$ where
$\pi_1$ and $\pi_2$ are the natural projections $\pi_1 : \mathbb{R}^n \rightarrow \mathbb{R}^k$ and $\pi_2 : \mathbb{R}^n
\rightarrow \mathbb{R}^{n-k}$, $\tilde{\omega}_1$ is a 1-form on a neighborhood of $0$ in $\mathbb{R}^k$ with an isolated
zero at the origin and $\tilde{\omega}_2$ is a radial 1-form on $\mathbb{R}^{n-k}$ at $0$.
\end{enumerate}

\begin{definition}
The radial index $\hbox{\em ind}_{Z,0}^\mathbb{R}\ \omega$ of the 1-form $\omega$ on $Z$ at $0$ is the sum:
$$1+\sum_{i=0}^r \sum_{Q \vert  \tilde{\omega}_{\vert S_i}(Q)=0} \hbox{\em ind}_{PH}(\tilde{\omega},Q,S_i),$$
where $\hbox{\em ind}_{PH}(\tilde{\omega},Q,S_i)$ is the Poincar\'e-Hopf index of the from $\tilde{\omega}_{\vert S_i}$ at $Q$ and where the sum is taken over all zeros of the 1-form $\tilde{\omega}$ on $(Z\setminus \{0\}) \cap B_\varepsilon$.
If $0$ is not a zero of $\omega$ on $Z$, we put $\hbox{\em ind}_{Z,0} ^\mathbb{R}\ \omega=0$.
\end{definition}
A straightforward corollary of this definition is that the radial index satisfies the law of conservation of number (see
Remark 9.4.6 in \cite{BSS} or the remark before Proposition 1 in \cite{EG1}).





Let us go back to the complex case. As in Section 2, $(X,0) \subset (\mathbb{C}^N,0)$ is an equidimensional reduced complex analytic germ of dimension $d$ in an open set $U \subset \mathbb{C}^N$. Let $\omega$ be a complex $1$-form on $U$ with an isolated singular point on $X$ at the origin.
\begin{definition}
The complex radial index ${\rm ind}_{X,0}^\mathbb{C}\ \omega$ of the complex $1$-form $\omega$ on $X$ at the origin is $(-1)^{d}$ times the index of the real $1$-form given by the real part of $\omega$.
\end{definition}

Let us denote $n_{i}=(-1)^{d-d_i-1} \left(\chi \big(\hbox{lk}^\mathbb{C}(V_i,X) \big)-1 \right)$, where $\{V_{i}\}$ is the Whitney stratification of $(X,0)$ considered in Section 2. In particular for a open stratum $V_{i}$ of $X$, $\hbox{lk}^\mathbb{C}(V_i,X)$ is empty and so $n_{i}=1$.
Let us define the Euler obstruction ${\rm Eu}_{Y,0}(\omega)$ to be equal to $1$ for a zero-dimensional variety $Y$.
Under this conditions Ebeling and Gusein-Zade proved in \cite{EG1} the following result which relates the radial index of a $1$-form to Euler obstructions.

\begin{theorem}\label{Th4} Let $(X,0) \subset (\mathbb{C}^{N},0)$ be the germ of a reduced complex analytic space at the origin, with a Whitney stratification $\{V_{i}\}$, $i= 0, \dots, q$, where $V_{0}=\{0\}$ and $V_{q}$ is the regular part of $X$. Then
$${\rm ind}_{X,0}^\mathbb{C}\ \omega=\sum_{i=0}^{q} n_{i}\cdot {\rm Eu}_{\overline{V_{i}},0}\ \omega.$$
\end{theorem}

\section{A lemma in subanalytic geometry}

Let $Y \subset \mathbb{R}^n$ be a closed subanalytic set equipped with a locally finite subanalytic Whitney stratification $\{S_i\}_{\alpha \in \Lambda}$: $Y= \cup S_i$.
Let $\rho : \mathbb{R}^n \rightarrow \mathbb{R}$ be a smooth subanalytic function such that $\rho ^{-1}(a)$ intersects $Y$ transversally (in the stratified sense). Then the following partition:
$$Y\cap \{ \rho  \le a \} = \bigcup  S_i \cap \{\rho <a \} \cup \bigcup S_i \cap \{\rho = a\},$$
is a Whitney stratification of the closed subanalytic set $Y \cap \{\rho  \le a \}$.

Let $\theta: \mathbb{R}^n \rightarrow \mathbb{R}$ be another smooth subanalytic function such that $\theta_{\vert Y \cap \{\rho \le a \} }$ admits an isolated critical point $p$ in $Y \cap \{\rho =a\}$ which is not a critical point of $\theta_{\vert Y}$.  If $S$ denotes the stratum of $Y$ that contains $p$, this implies that:
$$\nabla \theta_{\vert S}(p) =\lambda (p) \nabla \rho _{\vert S} (p),$$
with $\lambda(p) \not= 0$.
\begin{definition}
We will say that $p \in Y \cap \{ \rho  = a \}$ is an outward-pointing (resp. inward-pointing) critical point for $\theta_{\vert Y \cap \{ \rho  \le a \}}$ if $\lambda(p)>0$ (resp. $\lambda(p)<0$).
\end{definition}

Now let us suppose that $0 \in Y$ and that $\{0\}$ is a single stratum of $Y$. Let $\theta : U \subset \mathbb{R}^n \rightarrow \mathbb{R}$ be a smooth subanalytic function defined on an open neighborhood $U$ of $0$. We assume that $\theta_{\vert U \setminus \{0\} \cap S}$ has no critical point.

Let $\rho  :  \mathbb{R}^n \rightarrow \mathbb{R}$ be defined by $\rho (x)=x_1^2+\cdots+x_n^2$. It is known that for $\varepsilon >0$ small enough, the sphere $S_\varepsilon=\omega ^{-1}(\varepsilon^2)$ intersects $Y$ transversally. Let $p^\varepsilon$ be a critical point of $\theta_{\vert Y \cap S_\varepsilon}$. This means that there exists $\lambda(p^\varepsilon)$ such that:
$$ \nabla \theta_{\vert S}(p^\varepsilon) = \lambda(p^\varepsilon) \nabla \rho _{\vert S}(p^\varepsilon),$$
where $S$ is the stratum containing $p^\varepsilon$.
Note that $\lambda(p^\varepsilon) \not= 0$ since $\theta$ has no critical point on $U\setminus \{0\} \cap S$.

\begin{lemma}\label{CriticBoundary}
If $\varepsilon$ is small enough then $\theta(p^\varepsilon) \not= 0$. Furthermore $p^\varepsilon$ is outward-pointing (resp. inward-pointing) for $\theta_{\vert Y \cap B_\varepsilon}$ if and only if $\theta(p^\varepsilon)>0$
(resp. $\theta(p^\varepsilon)<0$).
\end{lemma}
{\it Proof.} If for $\varepsilon$ small enough, there is a critical point $p^\varepsilon$ of $\theta_{\vert Y \cap S^\varepsilon}$ such that $\theta(p^\varepsilon)=0$ then by the Curve Selection Lemma, there is a smooth subanalytic curve $p :[0,\nu[ \rightarrow Y$, $p(0)=0$, such that $p(t)$ is a critical point of $\theta_{\vert Y  \cap S_{\Vert p(t) \Vert}}$ and $g(p(t))=0$.
Since the stratification is locally finite, we can assume that $p(]0,\nu[)$ is included in a stratum $S$.
Hence, we have:
$$0=(\theta \circ p)'(t)=\langle \nabla \theta(p(t)), p'(t) \rangle =  \lambda(p(t)) \langle \nabla \rho _{\vert S}(p(t)),p'(t) \rangle,$$
because $p'(t)$ lies in $T_{p(t)} S$.  Therefore $(\rho \circ p)'(t)=0$ and $\rho\circ p$ is constant. But $\rho (p(t))$ tends to $0$ as $t$ tends to $0$ so $\rho \circ p$ is zero everywhere, which is a contradiction.

Now let us assume that $\lambda (p^\varepsilon) >0$. By the Curve Selection Lemma, there exists a smooth subanalytic curve $p:[0,\nu[ \rightarrow Y$ passing through $p^\varepsilon$ such that $p(0)=0$,
$p(]0,\nu[)$ is included in a stratum $S$ and for $t\not= 0$, $p(t)$ is a critical point of $\theta_{\vert S \cap S_{\Vert p(t) \Vert}}$ with $\lambda(p(t)) >0$. Therefore we have:
$$(\theta \circ p)'(t) = \lambda(p(t)) \langle \nabla \rho _{\vert S}(p(t)), p'(t) \rangle= \lambda(p(t)) (\rho  \circ p)'(t).$$
But $(\rho  \circ p)'>0$ for otherwise $(\rho  \circ p)' \le 0$ and $\rho \circ p$ would be decreasing.  Since $\rho (p(t))$ tends to $0$ as $t$ tends to $0$, this would imply that $\rho  \circ p(t) \le 0$, which is impossible. Hence we can conclude that $(\theta \circ p)'>0$ and $\theta \circ p$ is strictly increasing. Since $\theta \circ p(t)$ tends to $0$ as $t$ tends to $0$, we see that $\theta \circ p (t) >0$ for $t>0$. Similarly if $\lambda (p^\varepsilon) <0$ then $\theta(p^\varepsilon)<0$.  $\hfill \Box$

\section{Critical points and topology of Milnor fibres}
In this section we show  L\^e-Greuel type formulae (Theorem \ref{teo3.2} and Corollary \ref{IntNumbGen}) and we give applications for sets with isolated singularity.

As in section 2, $(X,0) \subset (\mathbb{C}^N,0)$ is an equidimensional reduced complex analytic germ of dimension $d$ in an open set $U$, equipped with a Whitney stratification $\{V_i\}$ such that $0$ belongs to the closure of all the strata. We write $X=\cup_{i=0}^q V_i$ where $V_0=\{0\}$ and $V_q=X_{\rm reg}$. We assume that the strata $V_0,\ldots,V_{q-1}$ are connected and that the analytic sets $\overline{V_0},\ldots,\overline{V_{q-1}}$ are reduced. We set $d_i={\rm dim}V_i$ for $i\in \{1,\ldots,q\}$.
Let $f : X \rightarrow \mathbb{C}$ be  a holomorphic function which is the restriction of a holomorphic function $F : U \rightarrow \mathbb{C}$. We assume that $f$ has an isolated singularity (or an isolated critical point) at $0$, i.e that $f$ has no critical point in a punctured neighborhood of $0$ in $X$. This implies that $X^f=X \cap f^{-1}(0)$ is a Whitney stratified set of dimension $d-1$, equipped with the stratification $\cup_{i=0}^q V_i \cap f^{-1}(0)= \cup_{i=0}^q V_i^f$.

Let us consider another holomorphic function $g : X \rightarrow \mathbb{C}$, restriction of a holomorphic function $G : U \rightarrow \mathbb{C}$. We also assume that $g$ has an isolated singularity at $0$ so that
$X^g=X\cap g^{-1}(0)= \cup_{i=0}^q V_i \cap g^{-1}(0)= \cup_{i=0}^q V_i^g$ is a Whitney stratification of $X^g$.
\begin{lemma}
The function $g : X^f \rightarrow \mathbb{C}$ has an isolated singularity at $0$ if and only if the function $f : X^g \rightarrow \mathbb{C}$ has an isolated singularity at $0$.
\end{lemma}
{\it Proof.} Let $\Sigma_{g_{\vert X^f}}$ denote the critical set of $g : X^f \rightarrow \mathbb{C}$. By the Curve Selection Lemma, it is easy to prove that $\Sigma_{g_{\vert X^f}}$ lies in $g^{-1}(0)$. Let $x$ be a point in $\Sigma_{g_{\vert X^f}}$ different from $0$. Then if $V(x)$ denotes the stratum of $X$ that contains $x$, we have:
$$dG_{\vert V(x)} (x) = \lambda(x) dF_{\vert V(x)}(x).$$
Since $x$ is not a critical point of $g$, $\lambda(x)$ is different from $0$ and:
$$dF_{\vert V(x)} (x) = \frac{1}{\lambda(x)} dG_{\vert V(x)}(x).$$
This last equality means that $x$ is a critical point of $f_{\vert X^g}$. $\hfill \Box$

For all $i \in \{1,\ldots,q \}$, we denote by $\Gamma^i_{f,g}$ the following relative polar set:
$$\Gamma^i_{f,g} = \left\{ x \in V_i \ \vert \ \hbox{rank} (dF_{\vert V_i}(x), dG_{\vert V_i}(x))< 2 \right\},$$
and we make the  assumption that $I_{X,0}(X^f,\overline{\Gamma^i_{f,g}}) < +\infty$, where $I_{X,0}(-,-)$ denotes the intersection multiplicity at the origin. This implies that the only critical point of $g_{\vert X^f}$ is $0$ and, by the previous lemma, that $f_{\vert X^g}$ has also an isolated singularity at $0$. Furthermore, it also implies that the number of critical points of $g_{\vert X \cap f^{-1}(\delta) \cap \mathring{B_\varepsilon}}$ is finite for $0 < \vert \delta \vert \ll \varepsilon \ll 1$.  Let us denote by $p_{i1},\ldots,p_{in_i}$ the critical points of $g_{\vert V_i \cap f^{-1}(\delta) \cap \mathring{B_\varepsilon} }$ and for each $j \in \{1,\ldots, n_i\}$, let us denote by $\mu_{ij}$ the Milnor number of $g_{\vert V_i}$ at $p_{ij}$. Our aim is to relate the $\mu_{ij}$'s to the topology of the Milnor fibres $X\cap f^{-1}(\delta) \cap B_\varepsilon$ and $X^g \cap f^{-1}(\delta) \cap B_\varepsilon$.
\begin{theorem}\label{teo3.2}
For $0 < \vert \delta \vert \ll \varepsilon \ll 1$, we have:
$$\chi \big( X \cap f^{-1}(\delta) \cap B_\varepsilon \big) -\chi \big( X^g \cap f^{-1}(\delta) \cap B_\varepsilon \big) =$$
$$ \sum_{i=1}^q (-1)^{d_i-1}  \sum_{j=1}^{n_i} \mu_{ij} \left( 1-\chi \big( \hbox{\em lk}^\mathbb{C}(V_i,X) \big) \right).$$
\end{theorem}
{\it Proof.}  Let us consider first the case when $X$ is a curve. In this situation, $X$ has two stratum $\{0\}$ and $X \setminus \{0\}$ and there is only one polar set $\Gamma_{f,g}^q$ which is exactly the stratum  $X \setminus \{0\}$. The critical points of $g_{\vert X \cap f^{-1}(\delta) \cap B_\varepsilon}$ are exactly the points in $X \cap f^{-1}(\delta) \cap B_\varepsilon$ and they have Milnor number $1$. Since $X^g \cap f^{-1}(\delta) \cap B_\varepsilon$ is empty as well as lk$^\mathbb{C} (X \setminus \{0\},X)$, the result is easy.

For the general case, we apply stratified Morse theory to the real and imaginary parts of $g$. Let us write $g=g_1 +\sqrt{-1}g_2$. Using the Cauchy-Riemann equations and local coordinates, it is not difficult to see that ${g_1}_{\vert X^f} : X^f \rightarrow \mathbb{R}$ and ${g_2}_{\vert X^f} : X^f \rightarrow \mathbb{R}$ have the same critical points as $g_{\vert X^f} : X^f \rightarrow \mathbb{C}$. Hence ${g_1}_{\vert X^f}$ and ${g_2}_{\vert X^f}$ have an isolated singularity at $0$. Similarly, ${g_1}_{\vert X \cap f^{-1}(\delta) \cap \mathring{B_\varepsilon}}$ and ${g_2}_{\vert X \cap f^{-1}(\delta) \cap \mathring{B_\varepsilon}}$ have the same critical points as
${g}_{\vert X \cap f^{-1}(\delta) \cap \mathring{B_\varepsilon}}$.

Let us study the critical points of $g_1$ on the stratified set $X \cap f^{-1}(\delta) \cap B_\varepsilon$. We can distinguish between two kinds of critical points : those lying in $X \cap f^{-1}(\delta) \cap \mathring{B_\varepsilon}$ and those lying in $X \cap f^{-1}(\delta) \cap S_\varepsilon$.  Applying Lemma \ref{CriticBoundary} to ${g_1}_{\vert X_f}$ and taking $\delta$ sufficiently close to $0$, we can say that the critical points of the second type satisfy the following properties:
\begin{enumerate}
\item they lie outside $\{ g_1 =0 \}$,
\item they are outward-pointing for ${g_1}_{\vert X \cap f^{-1}(\delta) \cap B_\varepsilon}$ in $\{g_1 > 0 \}$,
\item they are inward-pointing for ${g_1}_{\vert X \cap f^{-1}(\delta) \cap B_\varepsilon}$ in $\{g_1 < 0 \}$.
\end{enumerate}
Let $\tilde{g_1} : X \cap f^{-1}(\delta) \cap B_\varepsilon \rightarrow \mathbb{R}$ be a stratified Morse function close to $g_1$. Applying Theorem 2.3 of \cite{GMP} to the submersion $\mathbb{R}^N \times \mathbb{R}^N \rightarrow \mathbb{R}^N \times \mathbb{R}^N$, $(x,v) \mapsto G_1(x)+\sum_{i=1}^N x_iv_i$, where $G=G_1+\sqrt{-1}G_2$,  we can take $\tilde{g_1}$ to be the restriction of a real-analytic function.

For each $i \in \{1,\ldots, q\}$, $j \in \{1,\ldots,n_i\}$, let $\{q_{ij}^k \}$, $k \in \{1,\ldots,m_{ij} \}$, be the set of critical points of $\tilde{g_1}_{\vert V_i}$ lying close to $p_{ij}$. For each critical point $q_{ij}^k$, let $M_{\tilde{g_1}}(q_{ij}^k)$ be the local negative Milnor fibre of $\tilde{g_1}$ at $q_{ij}^k$. It is defined as follows:
$$ M_{\tilde{g_1}}(q_{ij}^k)= X \cap f^{-1}(\delta) \cap B_{\tilde{\varepsilon}}(q_{ij}^k) \cap \tilde{g_1}^{-1}\big(\tilde{g_1}(q_{ij}^k )- \nu \big),$$
where $0<\nu \ll \tilde{\varepsilon} \ll 1$ and $B_{\tilde{\varepsilon}}(q_{ij}^k)$ is the ball of radius $\tilde{\varepsilon}$ centered at $q_{ij}^k$.

Let $\tilde{\alpha}$ be a regular value of $\tilde{g_1}$ close to $0$. Applying Theorem 3.1 in \cite{Du}, we have:
$$\chi \big( X \cap f^{-1}(\delta) \cap B_\varepsilon \cap \{ \tilde{g_1} \ge \tilde{\alpha} \} \big) -\chi \big( X \cap f^{-1}(\delta) \cap B_\varepsilon \cap \{ \tilde{g_1}= \tilde{\alpha} \} \big)=$$
$$\sum_{i=1}^q \sum_{j=1}^{n_i} \sum_{k : \tilde{g_1}(q_{ij}^k) > \tilde{\alpha}} 1 -\chi \big( M_{\tilde{g_1}}(q_{ij}^k) \big).$$
Here we remark that the critical points of $\tilde{g_1}$ lying in $\{ \tilde{g_1} > \tilde{\alpha} \} \cap S_\varepsilon$ do not appear in the above equality. This is due to the fact that, since $\tilde{g_1}$ is close to $g_1$ and $\alpha$ close to $0$, these critical points are outward-pointing ${\tilde{g_1}}_{\vert X \cap f^{-1}(\delta) \cap B_\varepsilon}$. For such critical points, the local negative Milnor fibre has Euler characteristic $1$, as explained in \cite[Lemma 2.1]{Du}.

For each critical point $q_{ij}^k$, let $\lambda_{ij}^k$ be the Morse index of $\tilde{g_1}_{\vert V_i}$ at $q_{ij}^k$. Since $X \cap f^{-1}(\delta)$ is complex analytic, the normal Morse data of $\tilde{g_1}$ at $q_{ij}^k$ does not depend on $\tilde{g_1}$ nor on $q_{ij}^k$ and has the homotopy type of the pair $\left(\hbox{Cone}(\hbox{lk}^\mathbb{C}(V_i \cap f^{-1}(\delta) ,X\cap f^{-1}(\delta)),\hbox{lk}^\mathbb{C}(V_i\cap f^{-1}(\delta),X\cap f^{-1}(\delta)) \right)$ (see \cite[Corollary 1, p. 166]{GMP}). Moreover, since $\tilde{g_1}_{\vert X \cap f^{-1}(\delta) \cap B_\varepsilon}$ is a Morse function, the local Morse data at a critical point is the product of the tangential Morse data and the normal Morse data (see \cite[Section 3.7, p. 65]{GMP}). Hence, we can write:
$$ 1 -\chi \big( M_{\tilde{g_1}}(q_{ij}^k) \big)= (-1)^{\lambda_{ij}^k} \left(1-\chi \big(\hbox{lk}^\mathbb{C}(V_i\cap f^{-1}(\delta),X \cap f^{-1}(\delta) \big)\right).$$
Therefore, we have:
$$\chi \big( X \cap f^{-1}(\delta) \cap B_\varepsilon \cap \{ \tilde{g_1} \ge \tilde{\alpha} \} \big) -\chi \big( X \cap f^{-1}(\delta) \cap B_\varepsilon \cap \{ \tilde{g_1}= \tilde{\alpha} \} \big)=$$
$$\sum_{i=1}^q \sum_{j=1}^{n_i} \sum_{k : \tilde{g_1}(q_{ij}^k) > \tilde{\alpha}} (-1)^{\lambda_{ij}^k} \left(1-\chi \big(\hbox{lk}^\mathbb{C}(V_i\cap f^{-1}(\delta),X\cap f^{-1}(\delta)\big) \right).$$
Applying the same method to $-\tilde{g_1}$ and using the fact that the strata have even real dimension, we obtain:
$$\chi \big( X \cap f^{-1}(\delta) \cap B_\varepsilon \cap \{ \tilde{g_1} \le \tilde{\alpha} \} \big) -\chi \big( X \cap f^{-1}(\delta) \cap B_\varepsilon \cap \{ \tilde{g_1}= \tilde{\alpha} \} \big)=$$
$$\sum_{i=1}^q \sum_{j=1}^{n_i} \sum_{k : \tilde{g_1}(q_{ij}^k) < \tilde{\alpha}} (-1)^{\lambda_{ij}^k} \left(1-\chi \big(\hbox{lk}^\mathbb{C}(V_i\cap f^{-1}(\delta),X\cap f^{-1}(\delta) \big)\right).$$
Summing these two equalities and applying the Mayer-Vietoris sequence, we get:
$$\chi \big( X \cap f^{-1}(\delta) \cap B_\varepsilon  \big) -\chi \big( X \cap f^{-1}(\delta) \cap B_\varepsilon \cap \{ \tilde{g_1}= \tilde{\alpha} \} \big)=$$
$$\sum_{i=1}^q \sum_{j=1}^{n_i} \sum_{k =1}^{m_{ij}} (-1)^{\lambda_{ij}^k} \left(1-\chi \big(\hbox{lk}^\mathbb{C}(V_i\cap f^{-1}(\delta),X\cap f^{-1}(\delta)\big)\right).$$
But now for each $i \in \{1,\ldots,q\}$ and each $j \in \{1,\ldots,n_i\}$, $\sum_{k=1}^{m_{ij}} (-1)^{\lambda_{ij}^k} $ is the Poincar\'e-Hopf index of the form $d{g_1}_{\vert V_i \cap f^{-1}(\delta)}$. Since $g_1$ is the real part of $g$, this index is $(-1)^{d_i-1} \mu_{ij}$ (see for instance
 \cite[p. 235]{EG1}). Hence we have proved:
$$\chi \big( X \cap f^{-1}(\delta) \cap B_\varepsilon  \big) -\chi \big( X \cap f^{-1}(\delta) \cap B_\varepsilon \cap \{ \tilde{g_1}= \tilde{\alpha} \} \big)=$$
$$\sum_{i=1}^q \sum_{j=1}^{n_i}  (-1)^{d_i-1}  \mu_{ij}\left(1-\chi \big(\hbox{lk}^\mathbb{C}(V_i\cap f^{-1}(\delta),X\cap f^{-1}(\delta) \big) \right),$$
and therefore, since $\tilde{g_1}$ is close to $g_1$:
$$\chi \big( X \cap f^{-1}(\delta) \cap B_\varepsilon  \big) -\chi \big( X \cap f^{-1}(\delta) \cap B_\varepsilon \cap \{ g_1= \alpha \} \big)=$$
$$\sum_{i=1}^q \sum_{j=1}^{n_i}  (-1)^{d_i-1}  \mu_{ij} \left(1-\chi \big(\hbox{lk}^\mathbb{C}(V_i\cap f^{-1}(\delta),X\cap f^{-1}(\delta) \big) \right),$$
where $\alpha$ is a regular value of $g_1$ close to $0$. Applying Lemma \ref{complexlink} to $X$ and $f^{-1}(\delta)$, we obtain finally:
$$\chi \big( X \cap f^{-1}(\delta) \cap B_\varepsilon  \big) -\chi \big( X \cap f^{-1}(\delta) \cap B_\varepsilon \cap \{ g_1= \alpha \} \big)=$$
$$\sum_{i=1}^q \sum_{j=1}^{n_i}  (-1)^{d_i-1}  \mu_{ij} \left(1-\chi \big(\hbox{lk}^\mathbb{C}(V_i,X) \big)\right),$$

Now we are going to study the critical points of ${g_2}_{\vert X \cap f^{-1}(\delta) \cap B_\varepsilon \cap \{g_1 = \alpha \}}$. Using the Cauchy-Riemann equations and local coordinates, it is easy to see that ${g_2}_{\vert X \cap f^{-1}(\delta) \cap B_\varepsilon \cap \{g_1=\alpha\}}$ has no critical points in $\mathring{B_\varepsilon}$. Similarly ${g_2}_{\vert X^f \cap \{g_1=0\}}$ has an isolated singularity at the origin. Applying Lemma \ref{CriticBoundary} to ${g_2}_{\vert X^f \cap \{g_1=0\}}$ and taking $\delta$ and $\alpha$ very small, we can control the behaviour of the critical points of ${g_2}_{\vert X \cap f^{-1}(\delta) \cap B_\varepsilon \cap \{g_1=\alpha \}}$ lying in $S_\varepsilon$. Namely, we know that:
\begin{enumerate}
\item they lie outside $\{g_2=0 \}$,
\item they are outward-pointing for${g_2}_{\vert X \cap f^{-1}(\delta) \cap B_\varepsilon \cap \{g_1 = \alpha \}}$  in $\{g_2 >0 \}$,
\item they are inward-pointing for ${g_2}_{\vert X \cap f^{-1}(\delta) \cap B_\varepsilon \cap \{g_1 = \alpha \}}$ in $\{g_2 <0 \}$.
\end{enumerate}
Let $\beta$ be a small regular value of ${g_2}_{\vert X\cap f^{-1}(\delta) \cap \{g_1=\alpha \} \cap B_\varepsilon }$. Applying the same method as above, we find that:
$$\displaylines{
\qquad \chi \big( X \cap f^{-1}(\delta) \cap B_\varepsilon \cap \{ g_1= \alpha \} \cap \{ g_2 \ge \beta \}\big) - \hfill \cr
\hfill \chi \big( X \cap f^{-1}(\delta) \cap B_\varepsilon \cap \{ g_1= \alpha \} \cap \{ g_2= \beta \}\big)=0, \qquad \cr
}$$
$$\displaylines{
\qquad \chi \big( X \cap f^{-1}(\delta) \cap B_\varepsilon \cap \{ g_1= \alpha \} \cap \{ g_2 \le \beta \}\big) -\hfill \cr
\hfill \chi \big( X \cap f^{-1}(\delta) \cap B_\varepsilon \cap \{ g_1= \alpha \} \cap \{ g_2= \beta \}\big)=0. \qquad \cr
}$$
Hence, by the Mayer-Vietoris sequence, we get:
$$\chi \big( X \cap f^{-1}(\delta) \cap B_\varepsilon \cap \{ g_1= \alpha \}\big) =\chi \big( X \cap f^{-1}(\delta) \cap B_\varepsilon \cap \{ g_1= \alpha \} \cap \{ g_2= \beta \}\big).$$
Since $f_{\vert X^g}$ has an isolated singularity at the origin, $f^{-1}(\delta)$ intersects $X^g$ transversally. Hence, if $\alpha+\sqrt{-1}\beta$ is small enough, we have:
$$\chi \big( X \cap f^{-1}(\delta) \cap B_\varepsilon \cap \{ g= \alpha + \sqrt{-1}\beta \}\big) =
\chi \big( X \cap f^{-1}(\delta) \cap B_\varepsilon \cap \{ g=0\}\big)=$$
$$\chi \big( X^g \cap f^{-1}(\delta) \cap B_\varepsilon \big).$$
$\hfill \Box$

We can express this result in terms of intersection multiplicities.
\begin{corollary}\label{IntNumbGen}
For $0 < \vert \delta \vert \ll \varepsilon \ll 1$, we have:
$$\chi \big( X \cap f^{-1}(\delta) \cap B_\varepsilon \big) -\chi \big( X^g \cap f^{-1}(\delta) \cap B_\varepsilon \big) =$$
$$ \sum_{i=1}^q (-1)^{d_i-1} I_{X,0} (X^f,\overline{\Gamma^i_{f,g}})  \left( 1-\chi \big( \hbox{\em lk}^\mathbb{C}(V_i,X) \big) \right).$$
\end{corollary}
$\hfill \Box$

Let us apply these results to the case of a complex analytic variety with an isolated singularity. In this case, there are only two stratum $\{0\}$ and $X \setminus \{0\}$ and lk$^\mathbb{C}(X\setminus \{0\},X)$ is empty. Furthermore there is only one relative polar set:
$$\Gamma_{f,g} = \left\{ x \in X \setminus \{0\} \ \vert \ \hbox{rank} (dF_{\vert X \setminus \{0\}}(x), dG_{\vert X \setminus \{0\}}(x))< 2 \right\}.$$
\begin{corollary}\label{IntNumbIsol}
If $X$ has an isolated singularity at $0$ then for $0< \vert \delta \vert \ll \varepsilon  \ll 1$, we have:
$$\chi \big(X \cap f^{-1}(\delta) \cap B_\varepsilon \big)- \chi \big(X^g \cap f^{-1}(\delta) \cap B_\varepsilon \big)=(-1)^{d-1} I_{X,0}(X^f, \overline{\Gamma_{f,g}}).$$
\end{corollary}
$\hfill \Box$

\begin{corollary}
If $X$ is an ICIS defined by the equations  $f_1=\ldots=f_{N-d}=0$ such that $X\cap \{f=0\}$ and $X\cap \{f=0,g=0\}$ are also ICIS then the above equality gives that:
$$\mu(f_1,\ldots,f_{N-d},f)+\mu(f_1,\ldots,f_{N-d},f,g)=I_{X,0}(X^f,\overline{\Gamma_{f,g}}).$$
\end{corollary}

{\it Proof.}  We just use the equalities  $$\chi \big(X \cap f^{-1}(\delta) \cap B_\varepsilon \big)=1+(-1)^{d-1}\mu(f_1,\ldots,f_{N-d},f),$$ and $$\chi \big(X^g \cap f^{-1}(\delta) \cap B_\varepsilon \big)=1+(-1)^{d-2}\mu(f_1,\ldots,f_{N-d},f,g).$$ We recover a particular case of the L\^e-Greuel formula.
$\hfill \Box$

We can also apply our results to the case of a $1$-parameter smoothing. Namely let us assume that $X \subset \mathbb{C}^N$ has an isolated singularity at $0$ and let $\pi: \mathcal{X} \rightarrow D_\varepsilon$ be a flat morphism such that $\pi^{-1}(0)=X$ and $X_t=\pi^{-1}(t)$ is smooth for $t\not= 0$, where $\mathcal{X} \subset \mathbb{C}^M$ is a complex analytic set of dimension $d+1$ with an isolated singularity at $0$ and $D_\varepsilon$ is a small disk of radius $\varepsilon$ in $\mathbb{C}$. Let $g : \mathcal{X} \rightarrow \mathbb{C}$ be a holomorphic function, restriction of a holomorphic function $G : \mathbb{C}^M \rightarrow \mathbb{C}$. We assume that
$I_{\mathcal{X},0} (X, \overline{\Gamma_{\pi,g}}) < +\infty$ where $\Gamma_{\pi,g}$ is the relative polar set defined as above. With these conditions, Corollary \ref{IntNumbIsol} gives:
\begin{corollary}\label{Deformation}
For $0 < \vert t_0 \vert \ll \varepsilon \ll 1$, we have:
$$(-1)^d I_{\mathcal{X},0} (X,\overline{\Gamma_{\pi,g}})= \chi \big( \mathcal{X} \cap \pi^{-1}(t_0) \cap B_\varepsilon \big) -  \chi \big( \mathcal{X} \cap g^{-1}(0)\cap \pi^{-1}(t_0) \cap B_\varepsilon \big) .$$
\end{corollary}
$\hfill \Box$

When $X$ is a surface, $X \cap g^{-1}(0)$ is a curve and so:  $$\chi \big( \mathcal{X} \cap g^{-1}(0)\cap \pi^{-1}(t_0) \cap B_\varepsilon \big)=1-\mu(X \cap g^{-1}(0)).$$ If we define $\mu(X)$ to be the second Betti number of $X_{t_0} \cap B_\varepsilon$ then we obtain the following result:
$$ I_{\mathcal{X},0} ( X ,\overline{\Gamma_{\pi,g}}) +b_1(X_{t_0} \cap B_\varepsilon)=\mu(X)+\mu(X \cap g^{-1}(0)),$$
which was originally proved in \cite{PR}.


\section{L\^e-Greuel type formula and applications to curvatures}

In this section, we prove a L\^e-Greuel formula for the difference of the Euler obstruction and the Euler obstruction of a function, ${\rm B}_{f,X}(0)= {\rm Eu}_{X}(0) - {\rm Eu}_{f,X}(0)$, and we give applications to curvature integrals on the Milnor fibre.
We work with the objects and the assumptions of the previous section.

\begin{theorem} \label{IntNumbEuler}Assume that $f : X \rightarrow \mathbb{C}$  and $g :X \rightarrow \mathbb{C}$ have an isolated singularity at $0$ and that for all $i \in \{1,\ldots,q \}$, $I_{X,0}(X^f,\overline{\Gamma_{f,g}^i}) <+\infty$.
Then we have:
$$\displaylines{
\qquad {\rm B}_{f,X}(0)-{\rm B}_{f,X^{g}}(0) = (-1)^d I_{X,0} (X^f, \overline{\Gamma^q_{f,g}}).
 \qquad \cr
}$$
\end{theorem}

In order to prove the theorem, we need some lemmas relating Euler obstructions to complex links of strata.
\begin{lemma}
We have:
\begin{equation}{\rm Eu}_{X}(0)= 1 + \sum_{i=0}^{q-1}\left(\chi \big(\hbox{\em lk}^\mathbb{C}(V_i,X) \big)-1 \right){\rm Eu}_{\overline{V_{i}}}(0).\notag\end{equation}
\end{lemma}
{\it Proof.}
Let us consider $(x_{1}, x_{2}, \dots, x_{N})$ as complex coordinates of $\mathbb{C}^{N}$, where $x_{k}=u_{k}+\sqrt{-1} v_{k}$. This implies that $(u_{1},v_{1},\dots,u_{N},v_{N})$ are real coordinates of $\mathbb{R}^{2N}$. Let $\omega$ be a $1$-form defined by $\omega= \sum\overline{x_{k}}dx_{k}$, it means that:
 $$\omega= \sum (u_{k}-\sqrt{-1} v_{k}) (du_{k}+ \sqrt{-1} dv_{k}),$$
and so that:
 $$\omega = \sum (u_{k} du_{k} + v_{k}dv_{k}) + \sqrt{-1} \sum (u_{k} dv_{k}- v_{k} dv_{k}).$$

In this case, the real $1$-form ${\rm Re }\ \omega = \sum (u_{k}du_{k} + v_{k}dv_{k})$ is also a radial $1$-form, and  ${\rm ind}^{\mathbb{R}}_{X,0}\ {\rm Re}\ \omega =1$. Since ${\rm ind}_{X,0}^\mathbb{C}\ \omega = (-1)^{d}{\rm ind}^{\mathbb{R}}_{X,0} {\rm Re}\ \omega$, we find that:
$${\rm ind}_{X,0}^\mathbb{C}\ \omega = (-1)^{d}{\rm ind}^{\mathbb{R}}_{X,0} {\rm Re}\ \omega = (-1)^{d}.$$
As it was remarked before,
$${\rm Eu}_{X,0}\ \omega = (-1)^d {\rm Eu}_{X,0}{\rm Re}\ \omega .$$
Using this information and the definition of $n_{i}$ given in Section 3, we have the next equality:
$$n_{i} {\rm Eu}_{\overline{V_{i}},0}\ \omega= (-1)^{d-d_{i}-1}\left(\chi \big(\hbox{\rm lk}^\mathbb{C}(V_i,X) \big)-1\right)(-1)^{d_{i}}{\rm Eu}_{\overline{V_{i}}}(0).$$
Therefore, by Theorem \ref{Th4}  we conclude that:
$$(-1)^{d}=(-1)^{d}\left[ \sum_{i=0}^{q-1} \left(1-\chi \big(\hbox{\rm lk}^\mathbb{C}(V_i,X) \big)\right){\rm Eu}_{\overline{V_{i}}}(0)+{\rm Eu}_{X}(0) \right],$$
and so:
\begin{equation}{\rm Eu}_{X}(0)= 1 + \sum_{i=0}^{q-1}\left(\chi \big(\hbox{\rm lk}^\mathbb{C}(V_i,X) \big)-1\right){\rm Eu}_{\overline{V_{i}}}(0).\tag{1}\end{equation}
$\hfill \Box$

We have a similar result for the Euler obstruction of the function $f$.
\begin{lemma}
We have:
$$1-\chi(f^{-1}(\delta)\cap X \cap B_{\varepsilon})= \sum_{i=0}^{q}\left(1-\chi \big(\hbox{\em lk}^\mathbb{C}(V_i,X) \big) \right){\rm Eu}_{f,\overline{V_{i}}}(0).$$
\end{lemma}
{\it Proof}. On the one hand, applying the Theorem \ref{Th4} to the form $df$, we have:
$$\displaylines{
\qquad {\rm ind}_{X,0}^\mathbb{C} df = \sum_{i=0}^q n_{i}{\rm Eu}_{\overline{V_{i}},0}df =  \hfill \cr
\hfill (-1)^{d-d_i-1} \left(\chi \big(\hbox{\rm lk}^\mathbb{C}(V_i,X)\big)-1 \right) (-1)^{d_i}{\rm Eu}_{f,\overline{V_{i}}}(0). \qquad \cr
}$$
On the other hand, by Theorem 3 of \cite{EG1} we have:
$${\rm ind}_{X,0}^\mathbb{C} df=(-1)^{d}\left(1-\chi(f^{-1}(\delta)\cap X \cap B_{\varepsilon}) \right).$$
It follows that:
$$1-\chi(f^{-1}(\delta)\cap X \cap B_{\varepsilon})= \sum_{i=0}^{q}\left(1-\chi \big(\hbox{\rm lk}^\mathbb{C}(V_i,X) \big) \right){\rm Eu}_{f,\overline{V_{i}}}(0).$$
$\hfill \Box$

\begin{corollary}\label{FibMilnorEuler}
We have:
\begin{equation}\chi(f^{-1}(\delta)\cap X \cap B_{\varepsilon})=\sum_{i=0}^{q} \left(1-\chi \big(\hbox{\em lk}^\mathbb{C}(V_i,X) \big) \right) {\rm B}_{f,\overline{V_{i}}}(0) .\notag\end{equation}
\end{corollary}
{\it Proof.}
By the previous lemma, we  have the following equation:
\begin{equation}{\rm Eu}_{f,X}(0)=1-\chi(f^{-1}(\delta)\cap X \cap B_{\varepsilon})+ \sum_{i=0}^{q-1} \left(\chi \big(\hbox{\rm lk}^\mathbb{C}(V_i,X) \big)-1 \right){\rm Eu}_{f,\overline{V_{i}}}(0).\tag{2}\end{equation}
By the difference $(1)-(2)$ we arrive to:
$$\displaylines{
(3) \qquad {\rm B}_{f,X}(0)=\chi(f^{-1}(\delta) \cap X \cap B_{\varepsilon})+\hfill \cr
\hfill  \sum_{i=0}^{q-1}\left(\chi\big(\hbox{\rm lk}^\mathbb{C}(V_i,X)\big)-1  \right) {\rm B}_{f,\overline{V_{i}}},(0). \qquad \cr
}$$
Hence we find:
\begin{equation}\chi(f^{-1}(\delta)\cap X \cap B_{\varepsilon})=\sum_{i=0}^{q} \left(1-\chi \big(\hbox{\rm lk}^\mathbb{C}(V_i,X) \big) \right) {\rm B}_{f,\overline{V_{i}}}(0) .\notag\end{equation}
$\hfill \Box$

We are in position to prove the L\^e-Greuel formula for the Euler obstruction of $f$.

{\em Proof of Theorem \ref{IntNumbEuler}.} If $X$ is a curve, the result is easy because, by Theorem \ref{BMPS}:
$${\rm B}_{f,X}(0)=\chi(f^{-1}(\delta) \cap B_\varepsilon \cap X \setminus \{0\})= \chi(f^{-1}(\delta) \cap B_\varepsilon \cap X),$$
and ${\rm Eu}_{X^g}(0)={\rm Eu}_{f,X^g}(0)=1$. Thus it is enough to apply Corollary \ref{IntNumbIsol}.

Let us assume that ${\rm dim}\ X=d \ge 2$ and prove this result by induction on the depth of the stratification.
The first step is to consider the case when $X$ has isolated singularity at the origin. In this case our stratification will be $\{V_{0}=\{0\}, V_{1}=X_{\rm reg}\}$.
Applying Theorem \ref{BMPS}, we have:
$$\displaylines{
\qquad {\rm B}_{f,X}(0)-{\rm B}_{f,X^{g}}(0)= \hfill \cr
\hfill =\chi \big(X_{\rm reg}\cap B_{\varepsilon}\cap f^{-1}(\delta) \big) - \chi \big(X_{\rm reg}\cap B_{\varepsilon}\cap f^{-1}(\delta)\cap \{g=0\} \big)\qquad \cr
\hfill =\chi \big(X \cap B_{\varepsilon}\cap f^{-1}(\delta) \big) - \chi \big(X \cap B_{\varepsilon}\cap f^{-1}(\delta)\cap \{g=0\} \big).\qquad \cr
}$$
But, by Corollary \ref{IntNumbIsol}, we have:
$$\displaylines{
\qquad {\rm B}_{f,X}(0)-{\rm B}_{f,X^{g}}(0)= \hfill \cr
\hfill  (-1)^{d-1} I_{X,0} (X^f, \overline{\Gamma^1_{f,g}})\left(1-\chi \big({\rm lk}^{\mathbb{C}}(V_{1},X) \big) \right)= (-1)^{d-1} I_{X,0} (X^f, \overline{\Gamma^1_{f,g}}), \qquad \cr
}$$
because lk$^\mathbb{C}(V_1,X)$ is empty. Thus means that our assumption is true for the case of $X$ with isolated singularity at the origin.

Let us prove the general case. By the hypothesis of induction, for each $i \in \{1,\ldots,d-1\}$ we have:
$$\displaylines{
\qquad (-1)^{{d_{i}}-1}I_{X,0} (X^f, \overline{\Gamma^i_{f,g}})= {\rm B}_{f,\overline{V_{i}}}(0) - {\rm B}_{f,\overline{V_{i}}\cap\{g=0\}}(0). \qquad \cr
}$$
Using Corollary \ref{IntNumbGen}, we have:
$$ \sum_{i=1}^q (-1)^{d_i-1}  I_{X,0} (X^f, \overline{\Gamma^i_{f,g}})\left( 1-\chi \big( \hbox{\rm lk}^\mathbb{C}(V_i,X) \big) \right)=$$
$$\chi \big( X \cap f^{-1}(\delta) \cap B_\varepsilon \big) -\chi \big( X^g \cap f^{-1}(\delta) \cap B_\varepsilon \big).$$
Thus using the hypothesis of induction, we find:
$$\displaylines{
 \sum_{i=1}^{q-1} \left[{\rm B}_{f,\overline{V_{i}}}(0)- {\rm B}_{f,\overline{V_{i}}\cap\{g=0\}}(0)  \right]\left( 1-\chi \big( \hbox{\rm lk}^\mathbb{C}(V_i,X) \big) \right)+ \hfill \cr
\hfill (-1)^{d-1}I_{X,0} (X^f, \overline{\Gamma^q_{f,g}})=
\chi \big( X \cap f^{-1}(\delta) \cap B_\varepsilon \big) -\chi \big( X^g \cap f^{-1}(\delta) \cap B_\varepsilon \big). \qquad \cr
}$$
We can rewrite this equation as follows:
$$(-1)^{d-1}I_{X,0} (X^f, \overline{\Gamma^q_{f,g}})=   A- B,$$
where:
$$\displaylines{
\qquad A = \chi(f^{-1}(\delta)\cap X\cap B_{\varepsilon})- \sum_{i=0}^{q-1} {\rm B}_{f,\overline{V_{i}}} (0)  \left(1-\chi \big(\hbox{\rm lk}^\mathbb{C}(V_{i},X) \big) \right),  \hfill \cr
 \qquad \cr
}$$
and:
$$\displaylines{
\qquad B =\chi(f^{-1}(\delta)\cap X \cap B_{\varepsilon}\cap \{g=0\}) -  \hfill \cr
\hfill \sum_{i=0}^{q-1} {\rm B}_{f,\overline{V_{i}}\cap \{g=0\}}(0)  \left(1-\chi \big({\rm lk}^\mathbb{C}(V_{i}\cap \{g=0\},X \cap \{g=0\}) \big) \right), \qquad \cr
}$$
because by Lemma \ref{complexlink}, lk$^\mathbb{C}(V_i,X)={\rm lk}^\mathbb{C}(V_i \cap \{g=0\},X \cap \{g=0\})$.
Applying Corollary \ref{FibMilnorEuler}, we obtain:
\begin{equation}(-1)^{d-1}I_{X,0} (X^f, \overline{\Gamma^q_{f,g}})= {\rm B}_{f,X}(0)- {\rm B}_{f,X^g}(0). \notag\end{equation}
$\hfill \Box$


%

In \cite{LT1}, the authors show a formula to compute the Euler obstruction of $X$ at $0$ using polar multiplicities. Our next result generalizes this formula giving a similar formula for ${\rm B}_{f,X}(0)$, i.e., for the difference of the Euler obstruction and the Euler obstruction of the function. The strategy is to apply our results when $G$ is a generic linear function. Let $H \in \mathbb{C}P^{N-1}$ be a hyperplane defined by $H=\{ x \in \mathbb{C}^N \ \vert \ L(x)=0\}$ where $L$ is a linear function. It is well-known that if $H$ is general enough then $H$ intersects $X \setminus \{0\}$ and $X^f \setminus \{0\}$ transversally (see \cite[I.1.5.5]{Db2} for instance). Furthermore, if $l : X \rightarrow \mathbb{C}$ is the restriction to $X$ of $L$ then $\overline{\Gamma^q_{f,l}}$ is the general relative polar curve of the morphism $f : X \rightarrow \mathbb{C}$ (see \cite{LT1}  or \cite{HMS} for the definition of the general relative polar variety) and $I_{X,0}(X^f,\overline{\Gamma^q_{f,l}})<+\infty$ (see \cite{Lo}).
Following the notations of Loeser, we denote by $\Gamma_f^0$ the general relative polar curve of $f: X \rightarrow \mathbb{C}$. In this situation, we can apply Theorem \ref{IntNumbEuler}.
\begin{theorem}\label{TeissierLemmaSing}
We have:
$${\rm B}_{f,X}(0)-{\rm B}_{f,X\cap H}(0) = (-1)^{d-1} I_{X,0}(X^f,\Gamma_f^0),$$
where $H \in \mathbb{C}P^{N-1}$ is a generic hyperplane.
\end{theorem}
$\hfill \Box$

Since the left-hand side of the above equality and ${\rm Eu}_{X \cap H}$ do not depend on $H$ , we see that ${\rm Eu}_{f,X \cap H}(0)$ is independent on the choice of the generic hyperplane $H$. Furthermore, when $X=\mathbb{C}^N$, we recover the well-known Teissier lemma \cite{Te}, because in this situation: 
$${\rm B}_{f,X}(0)=1+(-1)^{N-1} \mu(f),$$
and
$${\rm B}_{f,X\cap H}(0)=1+(-1)^{N-2}  \mu'(f).$$


For $i \in \{0,\ldots,d-1\}$, let us denote by $\Gamma_f^i$ the general relative polar curve of the morphism $f : X \cap H^i \rightarrow \mathbb{C}$, where $H^i$ is a generic plane of codimension $i$ passing through $0$. The following theorem expresses ${\rm B}_{f,X}(0)$ in terms of the intersection multiplicities at the origin of the $\Gamma_f^i$'s with the $X \cap H^i$'s.
\begin{theorem}\label{EulerObstructionInterMult}
We have:
$${\rm B}_{f,X}(0)= \sum_{i=0}^{d-1} (-1)^{d-i-1} I_{X \cap H^i,0}(\Gamma_f^i, X^f \cap H^i).$$
\end{theorem}
{\it Proof.} Let $H^{d-1} \subset H^{d-2} \subset \cdots \subset H^1 \subset H^0=\mathbb{C}^N$ be a general flag such that $H^i$ has codimension $i$ for $i \in \{0,\ldots,d-1\}$. Denoting by $X^i$ the set $X \cap H^i$ and applying the previous theorem, we get:
$$\displaylines{
\qquad (-1)^{d-i-1} I_{X^i,0}(\Gamma_f^i, X^i \cap \{f=0 \})= {\rm B}_{f,X^i}(0)-{\rm B}_{f,X^{i+1}}(0), \qquad \cr
}$$
for $i \in \{0,\ldots,d-2 \}$. For $i=d-1$, $X^{d-1}$ is a curve and Theorem \ref{IntNumbEuler} gives:
$$I_{X^{d-1},0}(\Gamma_f^{d-1}, X^{d-1} \cap \{f=0\})= {\rm B}_{f,X^{d-1}}(0).$$
Summing all these equalities, we obtain the result. $\hfill \Box$

Applying this theorem when $f=l$ is the restriction to $X$ of a generic linear form leads to:
$${\rm Eu}_{X}(0) = \sum_{i=0}^{d-1} (-1)^{d-1-i} I_{X \cap H^i,0}(\Gamma_l^i, X^l \cap H^i).$$
But $I_{X \cap H^i,0}(\Gamma_l^i, X^l \cap H^i)$ is the multiplicity of the polar variety of codimension $d-i-1$ of $X$ (see \cite[Corollary 4.19]{LT1}). Hence we recover Corollary 5.12 of \cite{LT1}.
Based on this  last result, J.-P. Brasselet suggested that a possible definition for the Euler obstruction of a function could be the alternated sum of the relative polar multiplicities. This justifies why we choose to denote 
${\rm Eu}_{X}(0) - {\rm Eu}_{f,X}(0)$
by  ${\rm B}_{f,X}(0)$, which we call the Brasselet number. 




In \cite{K}, Kennedy presents a Gauss-Bonnet formula for the Milnor number of an analytic function $f : (\mathbb{C}^n,0) \rightarrow (\mathbb{C},0)$ with an isolated singularity. Using the integral formulas for the numbers $I_{X\cap H^i,0}(\Gamma_f^i, X^f \cap H^i)$ proved by Loeser \cite{Lo}, we can generalize Kennedy's formula to the singular case.   Let us recall first two objects used by Loeser.
For $i \in \{0,\ldots,d-1\}$, let $c^{w}_{d-1-i}(T_{f})$ be the ${(d-1-i)}$-th Chern Weil form associated with the relative tangent fiber bundle $T_{f}$ on the regular part of $X \cap f^{-1}(\delta)$ equipped with the hermitian structure given by the embedding of $X$ in $\mathbb{C}^{N}$.
Let us also define $\omega= \frac{\sqrt{-1}}{2\pi} \partial \overline{\partial} log \parallel z \parallel^{2}$ to be the inverse  image on $\mathbb{C}^{N}\setminus \{0\}$ of the K\"{a}hler form on $\mathbb{C}P^{N-1}$.
\begin{corollary}\label{KennedyGen}
We have the following integral formula:
$${\rm B}_{f,X}(0)=\lim_{\varepsilon \to 0} \lim_{\delta \to 0} \int_{X_{\rm reg} \cap f^{-1}(\delta) \cap B_{\varepsilon}} \sum_{i=0}^{d-1} c^{w}_{d-1-i}(T_{f})\wedge \omega^{i}.$$
\end{corollary}
{\it Proof.} By Theorem 1 in \cite{Lo}, we have that for $i \in \{0,\ldots,d-1\}$:
$$(-1)^{d-1-i}I_{X\cap H^i,0}(\Gamma_f^i, X^f \cap H^i) =\lim_{\varepsilon \to 0} \lim_{\delta \to 0} \int_{X_{\rm reg} \cap f^{-1}(\delta) \cap B_{\varepsilon}}  c^{w}_{d-1-i}(T_{f})\wedge \omega^{i}.$$
It is enough to apply the previous theorem. $\hfill \Box$

We can apply this last corollary to the stratum $V_{k}$, $k \in \{1,\ldots,q-1\}$, of $X$. Recalling that $d_k=\hbox{dim } V_k$, we have:
$$\lim_{\varepsilon \to 0} \lim_{\delta \to 0} \int_{V_{k} \cap f^{-1}(\delta) \cap B_{\varepsilon}} \sum_{i=0}^{d_k-1} c^{w}_{d_k-1-i}(T_{f})\wedge \omega^{i}= {\rm B}_{f,\overline{V_{k}}}(0),$$
where we keep the notation $T_f$ for the relative tangent bundle on $V_k \cap f^{-1}(\delta)$.

Multiplying by $1-\chi({\rm lk}^\mathbb{C}(V_{k},X))$ and using Corollary \ref{FibMilnorEuler}, we obtain a Gauss-Bonnet type formula for $X \cap f^{-1}(\delta) \cap B_{\varepsilon}$.
\begin{corollary}\label{GaussBonnetGen} We have the following Gauss-Bonnet formula:
$$\displaylines{
\qquad \chi(f^{-1}(\delta)\cap X \cap B_{\varepsilon})=  \hfill \cr
\hfill \sum_{k=1}^{q} \left(1-\left(\chi \big(\hbox{\em lk}^\mathbb{C}(V_k,X) \big) \right) \right) \lim_{\varepsilon \to 0}\lim_{\delta \to 0}\int_{V_{k} \cap f^{-1}(\delta) \cap B_{\varepsilon}}\sum_{i=0}^{d_k-1}
c^{w}_{d_k-1-i}(T_{f})\wedge \omega^{i}. \qquad \cr
} $$
\end{corollary}
$\hfill \Box$

\end{document}